\documentclass[12pt, a4paper]{article}
\usepackage{graphicx}
\usepackage{epsf}
\usepackage{amssymb,amsmath,eufrak,amsthm,amscd}
\usepackage{graphics}
\usepackage{textcomp}
\usepackage{hyperref}
\usepackage{amscd}
\usepackage{enumerate}
\usepackage{soul}
\usepackage{tikz}
\usepackage{verbatim}
\usetikzlibrary{calc,through,intersections,backgrounds,arrows}
\makeatletter
\let\@fnsymbol\@arabic
\makeatother

\newtheorem{theorem}{Theorem}[section]
\newtheorem{lemma}[theorem]{Lemma}

\newtheorem*{statement*}{Statement}
\newtheorem*{theorem*}{Theorem}
\newtheorem*{lemma*}{Lemma}
\newtheorem*{fact*}{Fact}

\newtheorem*{HD-Hoffman-theorem*}{Theorem \ref{HD-Hoffman-theorem}}

\theoremstyle{definition}
\newtheorem{definition}[theorem]{Definition}
\newtheorem*{definition*}{Definition}

\newtheorem*{example*}{Example}

\newtheorem*{exercise*}{Exercise}
\newtheorem{proposition}[theorem]{Proposition}
\newtheorem*{proposition*}{Proposition}

\newtheorem*{corollary*}{Corollary}

\newtheorem*{claim*}{Claim}

\theoremstyle{remark}
\newtheorem{remark}[theorem]{Remark}
\newtheorem*{remark*}{Remark}

\makeatletter
\newcommand{\l@abcd}[2]{\hbox to\textwidth{#1\dotfill #2}}
\makeatother

\title{On the chromatic number \\ of a simplicial complex\footnote{This work is a part of the Ph.D. thesis being written at the Hebrew University of Jerusalem, Israel.}}

\author{Konstantin Golubev
  \thanks{E-mail: \texttt{kost.golubev@mail.huji.ac.il}}\\
  \textit{\normalsize{Department of Mathematics,}}\\
  \textit{\normalsize{The Hebrew University of Jerusalem}} \\
  \textit{\normalsize{Givat Ram, Jerusalem 91904, Israel.}}
  }

\begin{document}
\maketitle
\begin{abstract}
In \cite{Ho} A.J. Hoffman proved a lower bound on the chromatic number of a graph in the terms of the largest and the smallest eigenvalues of its adjacency matrix. In this paper, we prove a higher dimensional version of this result and give a lower bound on the chromatic number of a pure $d$-dimensional simplicial complex in the terms of the spectra of the higher Laplacian operators. 
\end{abstract}

\section{Introduction}

The chromatic number of a graph is the minimal number of colors needed to color its vertices in such a way that no edge is monochromatic. It has been studied extensively by different methods. One of them is the spectral method, i.e. bounding the chromatic number by means of the spectra of various operators defined on the graph. For example, an upper bound was given by H.S. Wilf, \cite{Wi}, and a lower bound by A.J. Hoffman, \cite{Ho}. See also~\cite{WE} by P. Wocjan and C. Elphick, which contains a generalization of the latter bound as well as an exposition of the known results in the area. An advantage of the spectral method is that the spectrum of an operator on a finite graph can be calculated in polynomial time, while the problem of finding the chromatic number of a graph is NP-complete.


One operator is the Laplacian of a graph. Given a graph $G$ the Laplacian $\Delta$ of $G$ is an operator on the space $C^{0}$ of the real-valued functions on the vertex set $V$ of $G$. On a function $f\in C^0$ it acts as
$$
\Delta f(v) = \deg{v}\cdot f(v) - \sum_{u\sim v} f(u) = \sum_{u\sim v} \left(f(v)-f(u)\right),
$$
where $v\in V$, $\deg{v}$ is the number of edges adjacent to $v$, and $u\sim v$ stands for the vertices $u$ and $v$ being connected by an edge in $G$.

The chromatic number $\chi(G)$ of the graph $G$ is connected to its independence number $\alpha(G)$. The independence number $\alpha(G)$ of the graph $G$ is the cardinality of the largest subset of vertices, such that no two vertices of this subset are connected by an edge. Then $\alpha(G)\chi(G)\geq n$, where $n$ is the cardinality of the vertex set of $G$. 

In 1970, A.J. Hoffman gave a lower bound on the chromatic number in the terms of the eigenvalues of the adjacency matrix of the graph. Here we give a version of it in the terms of the Laplacian.

\begin{theorem}\label{Hoffman-theorem}(A.J. Hoffman, 1970, \cite{Ho})
For a nonempty graph $G$ on $n$ vertices,
$$
\alpha(G) \leq \frac{\mu-k}{\mu}\cdot n, \text{ and hence, }\chi(G)\geq \frac{\mu}{\mu-k},
$$
where $\alpha(G)$ is the independence number of $G$, $\chi(G)$ is the chromatic number of $G$, $k$ is the minimal degree of a vertex of $G$ and $\mu$ is the largest eigenvalue of the Laplacian $\Delta$ on $G$.
\end{theorem}

\begin{proof}
Let $A\subset V$ be the largest independent subset of $V$, and $B = V\setminus A$ be its complement. Note that both $A$ and $B$ are nonempty, since the graph $G$ is nonempty. Consider the following function $f\in C^{0}$ on the vertex set $V$
$$
f(v)=\begin{cases}
   -|B|,&{\text{if } v\in A;}\\
    |A|,&{\text{if } v\in B.}
  \end{cases}
$$
Then
$$
\mu\geq \frac{\langle\Delta f,f\rangle}{\langle f,f\rangle} \geq 
 \frac{|A|\cdot k\cdot (|A|+|B|)^{2}}{|A|\cdot |B|^{2} + |B|\cdot |A|^{2}} = \frac{k\cdot n}{n-|A|},
$$
and hence
$$
|A|=\alpha(G)\leq \frac{\mu-k}{\mu}\cdot n.
$$
Since $\chi(G)\alpha(G)\geq n$, the second bound follows.
\end{proof}

An important application of the Hoffman bound is given in \cite{LPS}, where non-bi-partite Ramanujan graphs are constructed. These graphs are regular, i.e. all vertices have the same degree, and are shown there to have chromatic number of order $\sqrt{k}$, where $k$ is the degree of regularity, and girth greater or equal to $\frac{4}{3}\log_{k-1}{n}$, where $n$ is the number of vertices. Thus, an explicit construction of graphs of arbitrarily large girth and arbitrarily large chromatic number was given.

Hoffman's theorem also serves as a powerful tool in extremal combinatorics. For example, in~\cite{Lo} L. Lov\'asz reproved the Erd\H{o}s-Ko-Rado theorem on the maximum size of a uniform family of intersecting sets by showing that the Hoffman bound is sharp for the corresponding Kneser graph. Another example is a different proof of the Deza-Frankl theorem on the number of intersecting permutations of $n$ elements. In~\cite{DF}, by purely combinatorial considerations M. Deza and P. Frankl proved that the largest set of intersecting permutations is of size $(n-1)!$. The work~\cite{Re} of P. Renteln, where the largest eigenvalue of the Laplacian of the derangement graph was computed, implies via the Hoffman theorem a sharp bound on its independence number, thus reproving the Deza-Frankl theorem. 

In this paper we present a generalization of the Hoffman result to higher dimensions. As a generalization of the notion of a finite graph we take a finite pure $d$-dimensional abstract simplicial complex. That is a family $X$ of subsets a finite vertex set $V$ closed under taking subsets, such that every maximal subset in the family is of size $d+1$. The elements of $X$ are called \textit{faces}. The dimension of a face is its cardinality minus one. By the degree of a $j$-face we mean the number of $(j+1)$-faces containing it.

The chromatic number $\chi(X)$ of a simplicial complex $X$ is the least number of colors needed to color its vertices in such a way that no maximal face is monochromatic. This is also known as the \textit{weak} chromatic number of a complex. The independence number $\alpha(X)$ of a simplicial complex $X$ is the size of the largest subset of vertices such that no maximal face of the complex has all its vetices in this subset. As in the case of graphs, $\alpha(X)\chi(X)\geq |V|$, where $V$ is the vertex set of the complex.

The Laplacian operators for simplicial complexes were introduced by B.~Eckmann in \cite{Ec} generalizing the above definition for graphs. The reader is referred to Section~\ref{Notations-and-definitions} for precise definitions.

The main result of this paper is the following theorem.

\begin{theorem}\label{HD-Hoffman-theorem}
Let $X$ be a nonempty pure $d$-dimensional simplicial complex on a finite vertex set $V$ of size $n$. Then
$$
\alpha(X) \leq \frac{\mu_0\dots\mu_{d-1}-(k_{0}+1)(k_{1}+2)\dots (k_{d-2}+d-1)k_{d-1}}{\mu_0\dots\mu_{d-1}} \cdot n,
$$
and hence,
$$
\chi(X) \geq \frac{\mu_0\dots\mu_{d-1}}{\mu_0\dots\mu_{d-1}-(k_{0}+1)(k_{1}+2)\dots (k_{d-2}+d-1)k_{d-1}},
$$
where $\chi(X)$ is the chromatic number of $X$, and for $0\leq j\leq d-1$, $k_j$ is the minimal degree of a $j$-face of $X$, and $\mu_{j}$ is the largest eigenvalue of the $j$-th upper Laplacian operator $\Delta_{j}^{+}$ of $X$.
\end{theorem}

We note that in the one-dimensional case, when a simplicial complex is just a graph, the upper Laplacian operator $\Delta_{0}^{+}$ coincides with the Laplacian operator $\Delta$ of a graph defined earlier, and that Theorem~\ref{Hoffman-theorem} is a special case of Theorem \ref{HD-Hoffman-theorem}.

In Section \ref{Main-theorem}, we prove the main result. In Section \ref{Complexes-with-complete-$(d-1)$-skeleton}, we formulate a version of Theorem \ref{HD-Hoffman-theorem} for the case of a pure $d$-dimensional simplicial complex with a complete $(d-1)$-skeleton and give a different, shorter, proof.

\section{Notations and definitions}\label{Notations-and-definitions}

\textit{An abstract simplicial complex} $X$ on a finite vertex set $V$ is a collection of subsets of $V$ closed under taking subsets.
Elements of $X$ are called \textit{faces} of the complex.

\textit{The dimension} of a face $\tau\in X$ is its cardinality minus one, $\dim(\tau)=|\tau|-1$. A face of dimension $j$ is called a \textit{$j$-face}. The set of all $j$-faces of $X$ is denoted $X^{j}$, in particular, $X^{0}=V$. \textit{The dimension} of a complex, $\dim(X)$, is the largest dimension of its face.
A $d$-dimensional complex is called \textit{pure}, if all the maximal faces are of dimension $d$.

By $X(j)$ we denote the set of all faces of dimension at most $j$. It is a $j$-dimensional simplicial complex by itself and it is called the \textit{$j$-skeleton} of $X$.

\textit{The degree} of a $j$-face $\tau\in X^j$ is the number of $(j+1)$-faces containing it, i.e. $\deg(\tau) = |\{\sigma\in X^{j+1}\mid \tau\subset \sigma\}|$. 


For a subset $A\subseteq V$, we denote the set of $j$-faces with all its vertices in $A$ by
$$
X^{j}(A)=\left\{\{v_{0},\dots,v_{j}\}\in X^{j}\mid v_{l}\in A,\text{ for all } l=0,\dots,j\right\}.
$$

For a $j$-face $\tau\in X^j$ and disjoint subsets $A_1,\dots,A_m\subseteq V$ denote the set of all $(j+m)$-faces containing $\tau$ and having exactly one vertex in each of $A_1,\dots,A_m$ by
$$
X^{j+m}(\tau,A_1,\dots,A_m) = \left\{\tau\cup\{v_1,\dots,v_m\}\in X^{j+m}\mid \forall 1\leq i\leq m:\, v_i\in A_i\right\}.
$$

\paragraph{Chromatic and independence numbers}

A subset $A\subseteq V$ of vertices of a 
simplicial complex $X$ is called \textit{independent}, if there is no maximal face of $X$ with all vertices in $A$. The \textit{independence number} $\alpha(X)$ of a complex $X$ is the size of the largest independent set of vertices.

The \textit{chromatic number} $\chi(X)$ of a 
complex $X$ is the least integer $\chi$ such that vertices of $X$ can be partitioned into $\chi$ disjoint independent sets. Note that $|V|\leq \alpha(X)\cdot \chi(X)$.

In other words, the chromatic number of a simplicial complex $X$ is the least number of colors needed to color the \textit{vertices} of $X$ is such a way that no  maximal face of $X$ is monochoromatic, i.e. no maximal face has all its vertices colored in one color. It is also known as the weak chromatic number of $X$, while the chromatic number of the underlying graph $X(1)$ of $X$, denoted $\chi(X(1))$, is known as the strong chromatic number of $X$. Clearly, $\chi(X)\leq \chi(X(1))$. In fact, a stronger inequality holds.

\begin{proposition}
For a $d$-dimensional complex $X$
$$
\chi(X)\leq \left\lceil\frac{\chi(X(1))}{d}\right\rceil,
$$
where $\lceil \cdot \rceil$ denotes the ceiling function.
\end{proposition}
\begin{proof}
Let $\chi_{1} = \chi(X(1))$, $m = \left\lceil\frac{\chi_{1}(X)}{d}\right\rceil$, and $\chi_{d}=\chi(X)$. Let $V=A_{1}\sqcup\dots\sqcup A_{\chi_{1}}$ be a proper 1-coloring of $X$ in $\chi_{1}$ colors, i.e. there is no edge with both endpoints in $A_{j}$ for all $1\leq j\leq \chi_{1}$.

Let $B_{1} = A_{1}\sqcup\dots\sqcup A_{d}$, $B_{2} = A_{d+1}\sqcup\dots\sqcup A_{2d},\dots, B_{m} = A_{(m-1)\cdot d + 1}\sqcup \dots\sqcup A_{\chi_{1}}$.

Since for $1\leq j \leq \chi_1$, $A_j$ contains no $1$-face of $X$, $B_i$ contains no $d$-face of $X$ for all $1\leq i \leq m$, and therefore, $\chi_{d}\leq m$.
\end{proof}

\paragraph{Laplacian operators}
A $j$-face $A\in X^{j}$ with an ordering of its vertices is called \textit{oriented} and  denoted $[A]$. For $0\leq j \leq d$, denote by $C^j=C^{j}(X,\mathbb{R})$ the vector space of all real-valued antisymmetric functions on oriented $j$-faces of the complex $X$. That is, for a function $f\in C^{j}$, a $j$-face $F=\{v_{0},\dots,v_{j}\}\in X^{j}$ and a permutation $\pi\in S_{j+1}$ the following equality holds
$$
f([v_{\pi(0)},\dots,v_{\pi(j)}]) = sgn(\pi)f([v_{0},\dots,v_{j}]).
$$
An inner product on $C^{j}$ is defined as
$$
\langle f,g\rangle = \sum_{A\in X^{j}}f([A])g([A]).
$$
Note that the sum runs over the non-oriented $j$-faces of $X$, but the functions are evaluated on oriented ones. Since the functions $f$ and $g$ are anti-symmetric, an orientation of each face may be chosen arbitrarily, and $\langle f,g \rangle$ is well-defined.

For $0\leq j \leq d-1$, the \textit{$j$-th coboundary map} $\delta_{j}:C^{j}\to C^{j+1}$ is defined as
$$
(\delta_{j} f)([v_{0},\dots,v_{j+1}]) = \sum_{i=0}^{j+1} (-1)^{i}f([v_{0},\dots,\widehat{v_{i}},\dots,v_{j+1}]).
$$
Note that this defines a cochain complex, i.e. $\delta_{j+1}\circ \delta_{j} = 0$.

For $0\leq j \leq d-1$, the \textit{$j$-th boundary map} $\partial_{j}:C^{j+1}\to C^{j}$ is defined as
$$
(\partial_{j} f)([v_{0},\dots,v_{j}]) = \sum_{\substack{u\in V:\\ \{u,v_{0},\dots,v_{j}\}\in X^{j+1}}} f([u,v_{0},\dots,v_{j}]).
$$
The $j$-th boundary map is the adjoint of the $j$-th coboundary map.

The following operator on $C^{j}$ is called \textit{the upper $j$-th Laplacian} of the complex:
$$
\Delta_{j}^{+} = \partial_{j} \circ \delta_{j}.
$$

\begin{remark}
One can also consider the lower $j$-th Laplacian defined as $\Delta_{j}^{-} = \delta_{j-1}\circ \partial_{j-1}$, and the $j$-th Laplacian defined as the sum of the upper and the lower ones. Both lower and upper Laplacians are postitve semidefinite. The spectrum of $(j+1)$-th lower Laplacian coincides with the spectrum of $j$-th upper Laplacian up to the multiplicity of the zero eigenvalue. 

Note that in the case of graphs, the upper Laplacian operator on $C^{0}$ of a graph coincides with the Laplacian operator defined earlier in the introduction

For an exposition of theory of Laplacian operators of an abstract simplicial complex, we address the reader to \cite{HJ} and the references therein.
\end{remark}

In this paper we deal with the largest eigenvalue of the upper Laplacian operators. For $0\leq j \leq d-1$, we denote the largest eigenvalue of $\Delta^+_j(X)$ by $\mu_j=\mu_j(X)$.
In~\cite[Theorem 3.4]{HJ}, a lower bound for the largest eigenvalue of the Laplacian of a complex is provided. 

\begin{theorem}\label{Max-EV-lower-bound} (D. Horak, J. Jost, 2011, \cite{HJ}, Theorem 3.4)
Let $X$ be a pure $d$-dimensional 
simplicial complex. Then for each $0\leq j \leq d-1$,
$$
K_{j} + (j+1)\leq \mu_{j},
$$
where 
$K_j$ is the maximal degree of a $j$-face and $\mu_j$ is the largest eigenvalue of the $j$-th upper Laplacian operator $\Delta_{j}^{+}$ of $X$.
\end{theorem}

\section{Proof of the Main Theorem}\label{Main-theorem}
In this section we prove
\begin{HD-Hoffman-theorem*}
Let $X$ be a nonempty pure $d$-dimensional simplicial complex on a finite vertex set $V$. Then
$$
\alpha(X)  \leq \frac{\mu_0\dots\mu_{d-1}-(k_{0}+1)\dots (k_{d-2}+d-1)k_{d-1}}{\mu_0\dots \mu_{d-1}}\cdot n,
$$
and hence,
$$
\chi(X) \geq \frac{\mu_0\dots \mu_{d-1}}{\mu_0\dots\mu_{d-1}-(k_{0}+1)\dots (k_{d-2}+d-1)k_{d-1}},
$$
where for $0\leq j\leq d-1$, $k_j$ is the minimal degree of a $j$-face of $X$, and  $\mu_j$ is the largest eigenvalue of the $j$-th upper Laplacian operator $\Delta_{j}^{+}$ of $X$.
\end{HD-Hoffman-theorem*}
The proof of this theorem is mainly based on the following lemma.

\begin{lemma}\label{lemma-descend} 
Let $X$ be a nonempty pure $d$-dimensional simplicial complex on a finite vertex set $V$. For $0\leq i\leq d-1$, denote by $k_i$ the minimal degree of an $i$-face of $X$. Let $0\leq j \leq d-2$, $A\subseteq V$ a subset of vertices, $B=V\setminus A$, and $0<\beta \leq k_{j+1}$ be a constant such that for each $(j+1)$-cell $\tau$ which has all its vertices in $A$, i.e. $\tau\in X^{j+1}(A)$, the following inequality holds
$$
|X^{j+2}(\tau,B)| \geq \beta.
$$ 
Then for each $j$-cell $\sigma\in X^{j}(A)$,
$$
|X^{j+1}(\sigma, B)| \geq \beta \cdot \frac{k_j+j+1}{\mu_{j+1}},
$$
where $\mu_{j+1}$ is the largest eigenvalue of the $(j+1)$-st upper Laplacian operator $\Delta_{j+1}^{+}$ of $X$.
\end{lemma}

\begin{proof}
Let $\sigma=\{v_0,v_1,\dots,v_j\}\in X^j(A)$ be a $j$-face of $X$ with all vertices in $A$. Denote by $a_{\sigma}$ the number of $(j+1)$-cells containing $\sigma$ and having the remaining vertex in $A$, and by $b_{\sigma}$ the number of $(j+1)$-cells containing $\sigma$ and having the remaining vertex in the complement $B=V\setminus A$, i.e.
$$
a_{\sigma} = |X^{j+1}(\sigma,A)|,\text{ and } b_{\sigma} = |X^{j+1}(\sigma,B)|. 
$$
Then
$$
a_{\sigma} + b_{\sigma} \geq k_j.
$$

Note that $b_{\sigma}> 0$, since $\sigma$ is contained in some $(j+1)$-face and $\beta >0$. Assume $a_{\sigma}=0$. Note that $k_{d-1}<k_{d-2}<\dots<k_0$, since $X$ is $d$-dimensional and pure, hence $\beta<k_j$. Also, $ \beta \leq K_{j+1}$, where $K_{j+1}$ is the maximal degree of a $(j+1)$-face of $X$. Then, by Theorem~\ref{Max-EV-lower-bound},
$$
\mu_{j+1} \geq K_{j+1} + (j+2) \geq \beta + (j+1) \frac{\beta}{k_{j}}
$$
and hence,
$$
b_{\sigma} \geq k_j \geq \beta \frac{k_{j}+j+1}{\mu_{j+1}}.
$$

Now assume that both $a_\sigma \neq 0$ and $b_\sigma \neq 0$. Denote the vertices in $A$ which form a $(j+1)$-cell with $\sigma$ by $u_1,\dots,u_{a_{\sigma}}$, and those in $B$ which form a $(j+1)$-cell with $\sigma$ by $w_1,\dots,w_{b_{\sigma}}$. Define a function $f\in C^{j+1}(X)$ in the following way:
\begin{itemize}
\item[--] For $1\leq i \leq a_\sigma$, 
$$
f[v_0,\dots, v_j,u_i]= (-1)^{j+2} \cdot b_{\sigma}.
$$
\item[--] For $1\leq i \leq b_\sigma$, 
$$
f[v_0,\dots, v_j,w_i]= (-1)^{j+1} \cdot a_{\sigma}.
$$
\item[--] For $1\leq i_1 \leq a_\sigma$, $1\leq i_2 \leq b_\sigma$ and $0\leq r \leq j$, 
$$
f[v_0,\dots,\widehat{v_r},\dots,v_j, u_{i_1},w_{i_2}]= (-1)^r.
$$
\item[--] For all other $(j+1)$-cells $\tau\in X^{j+1}$, $f[\tau]=0$.
\end{itemize}
Then for every $(j+2)$-cell of the form $\{v_0,\dots,v_j,u_{i_1},w_{i_2}\}$ the following equalities hold
$$
(\delta_{j+1} f)[v_0,\dots,v_j,u_{i_1},w_{i_2}]
$$
$$
 = \sum_{r=0}^j (-1)^r f[v_0,\dots,\widehat{v_r},\dots,v_j, u_{i_1},w_{i_2}]
+ (-1)^{j+1}f[v_0,\dots,v_j,w_{i_1}] 
$$
$$
+ (-1)^{j+2}f[v_0,\dots,v_j,u_{i_1}] 
$$
$$
 = \sum_{r=0}^j (-1)^{2r} + (-1)^{2(j+1)}a_{\sigma} + (-1)^{2(j+2)}b_{\sigma}= j+1 + a_{\sigma} + b_{\sigma}. 
$$
There are $a_{\sigma}$ $(j+1)$-cells in $A$ containing $\sigma$, hence there are at least $a_\sigma \beta$ $(j+2)$-cells of the form $\{v_0,\dots,v_j,u_{i_1},w_{i_2}\}$, i.e.,
$$
|X^{j+2}(\sigma,A,B)|\geq a_{\sigma} \beta.
$$
Therefore,
$$
\langle\Delta_{j+1} f,f \rangle = \langle \delta_{j+1} f, \delta_{j+1} f \rangle \geq a_\sigma \beta (j+1 + a_{\sigma} + b_{\sigma})^2.
$$ 

There are exactly $a_\sigma$ faces of the form $\{v_0,\dots, v_j, u_i\}$, and exactly $b_\sigma$ faces of the form $\{v_0,\dots,v_j,w_i\}$, and for each $0\leq r\leq j$ there are not more than $a_\sigma b_\sigma$ faces of the form $\{v_0,\dots,v_{r-1},v_{r},\dots,v_j, u_{i_1},w_{i_2}\}$. Hence
$$
\langle f,f \rangle \leq a_\sigma b_\sigma^2 + b_\sigma a_\sigma^2 + (j+1) a_\sigma b_\sigma = a_\sigma b_\sigma (j+1 + a_{\sigma} + b_{\sigma})
$$
Since $\mu_{j+1}$ is the largest eigenvalue of $\Delta_{j+1}$, we get that
$$
\mu_{j+1} \geq \frac{\langle\Delta_{j+1} f,f \rangle}{\langle f,f \rangle } 
\geq \frac{a_\sigma \beta (j+1 + a_{\sigma} + b_{\sigma})^2}{a_\sigma b_\sigma (j+1 + a_{\sigma} + b_{\sigma})} 
= \frac{\beta (j+1 + a_{\sigma} + b_{\sigma})}{b_\sigma}
$$
$$
\geq \frac{\beta (j+1 + k_j)}{b_\sigma},
$$
and hence
$$
b_\sigma = |X^{j+1}(\sigma,B)| \geq \beta\frac{k_j+j+1}{\mu_{j+1}}.
$$

\end{proof}
Now we turn to the proof of the Theorem~\ref{HD-Hoffman-theorem}.
\begin{proof}
Let $A\subseteq V$ be the largest independent set of vertices of $X$, in particular, there is no $d$-face with all vertices in $A$. Denote the complement $B=V\setminus A$ and $a=|A|$, $b=|B|$. Note that both $a$ and $b$ are positive. The desired inequality on the independence number is obtained by bounding from above and below the number $|X^{1}(A,B)|$ of 1-faces, i.e. edges, with exactly one vertex in $A$.

\paragraph{Upper bound:} Consider the function $f_{0}\in C^{0}$ on the vertices of $X$:
$$
f_{0}(v)=\begin{cases}
   -b,&{\text{if } v\in A;}\\
    a,&{\text{if } v\in B.}
  \end{cases}
$$
Since $\mu_0$ is largest eigenvalue of the $0$-th Laplacian $\Delta_0$,
$$
\mu_0 \geq \frac{\langle\Delta_{0} f_{0},f_{0}\rangle}{\langle f_{0},f_{0}\rangle} = \frac{\langle\delta_{0}f_{0},\delta_{0}f_{0}\rangle}{\langle f_{0},f_{0}\rangle}
$$
$$
= \frac{|X^{1}(A,B)|\cdot (a+b)^{2}}{a\cdot b^{2} + b\cdot a^{2}} = \frac{|X^{1}(A,B)|\cdot n}{(n-a) a },
$$
and hence
\begin{equation}\label{main-theorem-upper-bound}
\mu_0 \cdot \frac{(n-a)\cdot a}{n} \geq |X^{1}(A,B)|.
\end{equation}

\paragraph{Lower bound:} Since $A$ is independent, for each $(d-1)$-face in $A$, $\tau\in X^{d-1}(A)$ the following inequality holds
$$
|X^d(\tau,B)| \geq k_{d-1}.
$$
By a consecutive application of Lemma~\ref{lemma-descend}, we have that for each $0$-face, i.e. a vertex $v\in A$,
$$
|X^1(v,B)| \geq k_{d-1} \cdot \frac{k_{d-2} + d-1}{\mu_{d-1}} \dots \frac{k_0 + 1}{\mu_1}. 
$$ 
Hence,
\begin{equation}\label{main-theorem-lower-bound}
|X^1(A,B)| = \sum_{v\in A} |X^1(v,B)| \geq \frac{(k_0+1)\dots(k_{d-2}+d-1) k_{d-1}}{\mu_1\dots \mu_{d-1}} \cdot a.
\end{equation}

\paragraph{Conclusion:} By combining the bounds~\eqref{main-theorem-upper-bound} and \eqref{main-theorem-lower-bound}, we have 
$$
\mu_0 \cdot \frac{(n-a)\cdot a}{n} \geq |X^1(A,B)| \geq \frac{(k_0+1)\dots(k_{d-2}+d-1) k_{d-1}}{\mu_1\dots \mu_{d-1}} \cdot a
$$
which implies,
$$
a=\alpha(X)  \leq \frac{\mu_0\dots\mu_{d-1}-(k_{0}+1)\dots (k_{d-2}+d-1)k_{d-1}}{\mu_0\dots \mu_{d-1}}\cdot n.
$$
And since $\chi(X) \alpha(X) \geq n$, the second bound of the theorem follows.
\end{proof}

\section{Complexes with complete $(d-1)$-skeleton}\label{Complexes-with-complete-$(d-1)$-skeleton}

Denote $K_n^d$ the complete $d$-dimensional complex on $n$ vertices, i.e. every subset of size $(d+1)$ forms a face in $K_n^d$. A $d$-dimensional simplicial complex $X$ on a vertex set $V$, with $|V|=n$, is said to have a complete $(d-1)$-skeleton, if $X(d-1)=K_n^{d-1}$. 

\begin{theorem}\label{Complete-skeleton-Hoffman}
Let $X$ be a nonempty pure $d$-dimensional 
simplicial complex with complete $(d-1)$-skeleton on a finite vertex set $V$ of size $n$, then
$$
\alpha(X) \leq \frac{\mu-k_{d-1}}{\mu}\cdot n,
$$
and
$$
\chi(X) \geq \frac{\mu}{\mu-k_{d-1}},
$$
where $k_{d-1}$ is the minimal degree of a $(d-1)$-face, $\mu=\mu_{d-1}$ is the largest eigenvalue of the $(d-1)$-th upper Laplacian operator $\Delta_{d-1}^{+}$ of $X$.
\end{theorem}
Since for each $j=0,\dots,d-2$ all $j$-faces are of the same degree $n-(j+1)$, and the largest eigenvalue of $\Delta_{j}^{+}$ is
$$
\mu_{j} = n,\,\, 0\leq j \leq d-2,
$$
the proposition follows directly from the main theorem. Nevertheless, this special case can be proved directly and differently as follows.

\begin{proof}
As before, let $A\subseteq V$ be the largest independent set of vertices of $X$, let $B = V\setminus A$ be its complement, and $a=|A|$, $b=|B|=n-a$. Note that, $a\geq d$. Partition $A$ into a disjoint union of $d$ non-empty subsets $A=A_{0}\sqcup\dots\sqcup A_{d-1}$, and denote $|A_{j}| = a_{j}, j=0,\dots,d-1.$ Consider the following function $f\in C^{d-1}(X,\mathbb{R})$
\begin{itemize}
\item[--] For $v_i\in A_i,\,\, 0\leq i \leq d-1$, 
$$
f[v_0,\dots, v_{d-1}]= b.
$$
\item[--] For $v_i\in A_i,\,\, 0\leq i \leq d-1$, $\,u\in B$ and $0\leq l \leq d-1$,
$$
f[v_0,\dots,v_{l-1},v_{l+1},\dots, v_{d-1},u]= (-1)^l \cdot a_{l}.
$$
\item[--] For all other $(d-1)$-cells $\tau\in X^{j+1}$, $f[\tau]=0$.
\end{itemize}
Then
$$
\langle f,f \rangle = a_{0}\cdots a_{d-1}\cdot b \cdot (b+a_{0}+\dots a_{d-1}) = a_{0}\cdots a_{d-1}\cdot (n-a) \cdot n
$$
and, since $A$ is independent and $X$ has complete skeleton,
$$
\langle \Delta_{d-1} f,f \rangle 
\geq a_{0} \cdots a_{d-1} \cdot k_{d-1} \cdot (b+a_{0}+\dots a_{d-1})^{2} 
= a_{0} \cdots a_{d-1} \cdot k_{d-1} \cdot n^{2}
$$
And since
$$
\mu_{d-1} \geq \frac{\langle\Delta_{d-1} f,f\rangle}{\langle f,f\rangle} = \frac{k_{d-1}\cdot n}{n-a}
$$
the bound follows.
\end{proof}

This proof was inspired by the work~\cite{PRT} of O. Parzanchevski, R. Rosenthal and R. Tessler on a generalization of the Cheeger inequality.
\smallskip

\begin{remark}
The bound presented in Theorem~\ref{Complete-skeleton-Hoffman} is sharp for the independence number of the full $d$-dimensional complex $K_{n}^{d}$ on $n$ vertices.

The complex $K_{n}^{d}$ has a full $(d-1)$-skeleton and every $(d-1)$-face is contained in $k_{d-1} = n-d$ $d$-faces. The largest eigenvalue of $\Delta_{d-1}^{+}$ is $\mu_{d-1} = n$.
The independence number of $K_{n}^{d}$ is equal to $d$, and it is equal to the bound given by Theorem \ref{HD-Hoffman-theorem}
$$
\alpha_{d}(K_{n}^{d}) = d = \frac{\mu_{d-1}-k_{d-1}}{\mu_{d-1}}\cdot n.
$$
\end{remark}

\paragraph{Acknowledgements.}
The author is grateful to his advisor A. Lubotzky as well as to Oren Becker, Shai Evra, Ori Parzanchevski, Ron Rosenthal and Benny Sudakov for fruitful discussions on the topic, and to ERC for the support. 

The author is grateful to two anonymous referees, whose valuable remarks helped shortening the proofs and enriching the content of the paper.

The final draft of the paper was prepared during a short visit of the author to ITC-ETHZ.

\end{document}